\documentclass{article}
\usepackage{amsfonts}
\usepackage{amsmath}
\usepackage{harvard}

\newtheorem{theorem}{Theorem}

\newtheorem{corollary}[theorem]{Corollary}

\newtheorem{definition}[theorem]{Definition}
\newtheorem{example}[theorem]{Example}

\newtheorem{remark}[theorem]{Remark}

\newenvironment{proof}[1][Proof]{\noindent\textbf{#1.} }{\ \rule{0.5em}{0.5em}}
\numberwithin{theorem}{section}
\numberwithin{equation}{section}
\input{tcilatex}

\begin{document}

\title{Jet Geometrical Objects Produced by Linear ODEs Systems and Superior
Order ODEs}
\author{Mircea Neagu \\
{\small July 2007; Revised June 2009 (minor revisions)}}
\date{}
\maketitle

\begin{abstract}
The aim of this paper is to construct a Riemann-Lagrange geometry on 1-jet
spaces, in the sense of d-connections, d-torsions, d-curvatures,
electromagnetic d-field and geometric electromagnetic Yang-Mills energy,
starting from a given linear ODEs system or a given superior order ODE. The
case of a non-homogenous linear ODE of superior order is disscused.
\end{abstract}



\section{Introduction}

\hspace{4mm} According to Olver's opinion expressed in [7] and in private
discussions, we point out that the 1-jet spaces are main mathematical models
necessary for the study of classical or quantum field theories. In a such
context, the contravariant differential geometry of the 1-jet spaces was
intensively studied by authors like Asanov [1] or Saunders [9].

Situated in the direction initiated by Asanov [1], it has been recently
developed the \textit{Riemann-Lagrange geometry of 1-jet spaces} [2], [4],
which is a geometrical theory on 1-jet spaces analogous with the well known 
\textit{Lagrange geometry of the tangent bundle} developed by Miron and
Anastasiei [3].

It is important to note that the Riemann-Lagrange geometry of the 1-jet
spaces allows the regarding of the solutions of a given ODEs (respectively,
PDEs) system as \textit{geodesics} [10] (respectively, \textit{generalized
harmonic maps} [6] or \textit{potential maps} [11]) in a convenient
Riemann-Lagrange geometrical structure on 1-jet spaces. In this way, it was
given a final solution for an open problem suggested by Poincar\'{e} [8] (%
\textit{find the geometric structure which transforms the field lines of a
given vector field into geodesics}) and generalized by Udri\c{s}te [10] (%
\textit{find the geometrical structure which converts the solutions of a
given first order PDEs system into harmonic maps}).

In this context, using the Riemann-Lagrange geometrical methods, it was
constructed an entire contravariant differential geometry on 1-jet spaces,
in the sense of d-connections, d-torsions, d-curvatures, electromagnetic
d-field and geometric electromagnetic Yang-Mills energy, starting only with
a given ODEs [5] (respectively, PDEs [6]) system of order one and a pair of
Riemannian metrics.

\section{Jet Riemann-Lagrange geometry produced by a non-linear ODEs system
of order one and a pair of Riemannian metrics}

\hspace{4mm} In this Section we present the main jet Riemann-Lagrange
geometrical ideas used for the geometrical study of a given non-linear first
order ODEs system. For more details, the reader is invited to consult the
works [4], [5] and [11].

Let $T=[a,b]\subset \mathbb{R}$ be a compact interval of the set of real
numbers and let us consider the jet fibre bundle of order one%
\begin{equation*}
J^{1}(T,\mathbb{R}^{n})\rightarrow T\times \mathbb{R}^{n},\text{ }n\geq 2,
\end{equation*}%
whose local coordinates $(t,x^{i},x_{1}^{i}),$ $i=\overline{1,n},$ transform
after the rules%
\begin{equation*}
\widetilde{t}=\widetilde{t}(t),\text{ }\widetilde{x}^{i}=\widetilde{x}%
^{i}(x^{j}),\text{ }\widetilde{x}_{1}^{i}=\frac{\partial \widetilde{x}^{i}}{%
\partial x^{j}}\frac{dt}{d\widetilde{t}}\cdot x_{1}^{j}.
\end{equation*}

\begin{remark}
From a physical point of view, in the 1-jet space of \textbf{physical events}
the coordinate $t$ has the physical meaning of \textbf{relativistic time},
the coordinates $(x^{i})_{i=\overline{1,n}}$ represent \textbf{spatial
coordinates} and the coordinates $(x_{1}^{i})_{i=\overline{1,n}}$ have the
physical meaning of \textbf{relativistic velocities}.
\end{remark}

Let $X=\left( X_{(1)}^{(i)}(t,x^{k})\right) $ be an arbitrary given d-tensor
field on the first order jet space $J^{1}(T,\mathbb{R}^{n})$, which produces
the jet non-linear ODEs system of order one (\textit{jet dynamical system})%
\begin{equation}
x_{1}^{i}=X_{(1)}^{(i)}(t,x^{k}(t)),\text{ }\forall \text{ }i=\overline{1,n},
\label{ODEs}
\end{equation}%
where $c(t)=(x^{i}(t))$ is an unknown curve on $\mathbb{R}^{n}$ and we use
the notations%
\begin{equation*}
x_{1}^{i}\overset{not}{=}\dot{x}^{i}=\frac{dx^{i}}{dt},\text{ }\forall \text{
}i=\overline{1,n}.
\end{equation*}

Suppose now that we fixed \textit{a priori }two Riemannian structures $%
(T,h_{11}(t))$ and $(\mathbb{R}^{n},\varphi _{ij}(x))$, where $x=(x^{k})_{k=%
\overline{1,n}}$, together with their attached Christoffel symbols $%
H_{11}^{1}(t)$ and $\gamma _{jk}^{i}(x)$. Automatically, the jet non-linear
ODEs system of order one (\ref{ODEs}), together with the pair of Riemannian
metrics%
\begin{equation*}
\mathcal{P}=(h_{11}(t),\varphi _{ij}(x)),
\end{equation*}%
produce the \textit{jet least squares Lagrangian function} 
\begin{equation*}
JLS_{\mathcal{P}}^{\text{ODEs}}:J^{1}(T,\mathbb{R}^{n})\rightarrow \mathbb{R}%
_{+},
\end{equation*}%
expressed by%
\begin{equation*}
JLS_{\mathcal{P}}^{\text{ODEs}}(t,x^{k},x_{1}^{k})=h^{11}(t)\varphi _{ij}(x)%
\left[ x_{1}^{i}-X_{(1)}^{(i)}(t,x)\right] \left[
x_{1}^{j}-X_{(1)}^{(j)}(t,x)\right] .
\end{equation*}

It is obvious that the \textit{global minimum points} of the \textit{jet
least squares energy action}%
\begin{equation*}
\mathbb{E}_{\mathcal{P}}^{\text{ODEs}}(c(t))=\int_{a}^{b}JLS_{\mathcal{P}}^{%
\text{ODEs}}(t,x^{k}(t),\dot{x}^{k}(t))\sqrt{h_{11}(t)}dt
\end{equation*}%
are exactly the solutions of class $C^{2}$ of the jet non-linear ODEs system
of order one (\ref{ODEs}). In other words, we have

\begin{theorem}
The solutions of class $C^{2}$ of the first order ODEs system (\ref{ODEs})
verify the second order Euler-Lagrange equations produced by the jet \textit{%
least squares Lagrangian function} $JLS_{\mathcal{P}}^{\text{ODEs}}$, namely
(\textbf{jet geometric dinamics})%
\begin{equation}
\frac{\partial \left[ JLS_{\mathcal{P}}^{\text{ODEs}}\right] }{\partial x^{i}%
}-\frac{d}{dt}\left( \frac{\partial \left[ JLS_{\mathcal{P}}^{\text{ODEs}}%
\right] }{\partial \dot{x}^{i}}\right) =0,\text{ }\forall \text{ }i=%
\overline{1,n}.  \label{E-L-P}
\end{equation}
\end{theorem}

\begin{remark}
Conversely, the above statement does not hold good because there exist
solutions for the second order Euler-Lagrange ODEs system (\ref{E-L-P})
which are not global minimum points for the jet least squares energy action $%
\mathbb{E}_{\mathcal{P}}^{\text{ODEs}}$, that is which are not solutions for
the jet first order ODEs system (\ref{ODEs}).
\end{remark}

As a conclusion, we believe that we may regard $JLS_{\mathcal{P}}^{\text{ODEs%
}}$ as a natural geometrical substitut on $J^{1}(T,\mathbb{R}^{n})$ for the
jet first order ODEs system (\ref{ODEs}).

But, we point out that a Riemann-Lagrange geometry on $J^{1}(T,\mathbb{R}%
^{n})$ produced by the jet least squares Lagrangian function $JLS_{\mathcal{P%
}}^{\text{ODEs}}$, via its second order Euler-Lagrange equations (\ref{E-L-P}%
), geometry in the sense of non-linear connection, generalized Cartan
connection, d-torsions and d-curvatures, is now completely done in the
papers [4], [5] and [6]. Moreover, a distinguished jet electromagnetic
2-form, characterized by some natural generalized Maxwell equations and a
geometric jet Yang-Mills energy [5], is constructed from the jet least
squares Lagrangian function $JLS_{\mathcal{P}}^{\text{ODEs}}$.

\begin{definition}
Any geometrical object on $J^{1}(T,\mathbb{R}^{n})$, which is produced by
the jet least squares Lagrangian function $JLS_{\mathcal{P}}^{\text{ODEs}}$,
via the Euler-Lagrange equations (\ref{E-L-P}), is called \textbf{%
geometrical object produced by the jet first order ODEs system (\ref{ODEs})
and the pair of Riemannian metrics }$\mathcal{P}$.
\end{definition}

In this context, we give the following jet Riemann-Lagrange geometrical
result, which is proved in [5] and, for the multi-time general case, in [6].
For more details, the reader is invited to consult the book [4].

\begin{theorem}
\label{MainThODEs} (i) The \textbf{canonical non-linear connection on }$%
J^{1}(T,\mathbb{R}^{n})$\textbf{\ produced by the jet first order ODEs
system (\ref{ODEs}) and the pair of Riemannian metrics} $\mathcal{P}$ is%
\begin{equation*}
\Gamma _{\mathcal{P}}^{\text{ODEs}}=\left(
M_{(1)1}^{(i)},N_{(1)j}^{(i)}\right) ,
\end{equation*}%
whose local components are given by%
\begin{equation*}
M_{(1)1}^{(i)}=-H_{11}^{1}x_{1}^{i}\text{ and }N_{(1)j}^{(i)}=\gamma
_{jk}^{i}x_{1}^{k}-\frac{1}{2}\left[ X_{(1)||j}^{(i)}-\varphi
^{ir}X_{(1)||r}^{(s)}\varphi _{sj}\right] ,
\end{equation*}%
where%
\begin{equation*}
X_{(1)||j}^{(i)}=\frac{\partial X_{(1)}^{(i)}}{\partial x^{j}}%
+X_{(1)}^{(m)}\gamma _{mj}^{i}.
\end{equation*}

(ii) The \textbf{canonical generalized Cartan connection }$C\Gamma _{%
\mathcal{P}}^{\text{ODEs}}$\textbf{\ produced by the jet first order ODEs
system (\ref{ODEs}) and the pair of Riemannian metrics} $\mathcal{P}$ has
the adapted components%
\begin{equation*}
C\Gamma _{\mathcal{P}}^{\text{ODEs}}=(H_{11}^{1},0,\gamma _{jk}^{i},0).
\end{equation*}

(iii) The effective adapted components of the \textbf{torsion} d-tensor 
\textbf{T}$_{\mathcal{P}}^{\text{ODEs}}$ of the canonical generalized Cartan
connection $C\Gamma _{\mathcal{P}}^{\text{ODEs}}$ \textbf{produced by the
jet first order ODEs system (\ref{ODEs}) and the pair of Riemannian metrics} 
$\mathcal{P}$ are%
\begin{equation*}
R_{(1)1j}^{(i)}=\frac{1}{2}\left[ X_{(1)||j//1}^{(i)}-\varphi
^{ir}X_{(1)||r//1}^{(s)}\varphi _{sj}\right]
\end{equation*}%
and%
\begin{equation*}
R_{(1)jk}^{(i)}=r_{jkm}^{i}x_{1}^{m}-\frac{1}{2}\left[ X_{(1)||j||k}^{(i)}-%
\varphi ^{ir}X_{(1)||r||k}^{(s)}\varphi _{sj}\right] ,
\end{equation*}%
where $r_{ijk}^{l}(x)$ are the components of the curvature tensor of the
Riemannian metric $\varphi _{ij}(x)$ and%
\begin{equation*}
\begin{array}{l}
X_{(1)||j//1}^{(i)}=\dfrac{\partial X_{(1)||j}^{(i)}}{\partial t}%
-X_{(1)||j}^{(i)}H_{11}^{1},\medskip \\ 
X_{(1)||j||k}^{(i)}=\dfrac{\partial X_{(1)||j}^{(i)}}{\partial x^{k}}%
+X_{(1)||j}^{(m)}\gamma _{mk}^{i}-X_{(1)||m}^{(i)}\gamma _{jk}^{m}.%
\end{array}%
\end{equation*}

(iv) The effective adapted components of the \textbf{curvature} d-tensor 
\textbf{R}$_{\mathcal{P}}^{\text{ODEs}}$ of the canonical generalized Cartan
connection $C\Gamma _{\mathcal{P}}^{\text{ODEs}}$ \textbf{produced by the
jet first order ODEs system (\ref{ODEs}) and the pair of Riemannian metrics} 
$\mathcal{P}$ are only $R_{ijk}^{l}=r_{ijk}^{l}.$

(v) The \textbf{geometric electromagnetic distinguished 2-form produced by
the jet first order ODEs system (\ref{ODEs}) and the pair of Riemannian
metrics} $\mathcal{P}$ has the expression%
\begin{equation*}
F_{\mathcal{P}}^{\text{ODEs}}=F_{(i)j}^{(1)}\delta x_{1}^{i}\wedge dx^{j},
\end{equation*}%
where%
\begin{equation*}
\delta x_{1}^{i}=dx_{1}^{i}+M_{(1)1}^{(i)}dt+N_{(1)k}^{(i)}dx^{k}
\end{equation*}%
and, if $h^{11}=1/h_{11}$, then%
\begin{equation*}
F_{(i)j}^{(1)}=\frac{h^{11}}{2}\left[ \varphi _{im}X_{(1)||j}^{(m)}-\varphi
_{jm}X_{(1)||i}^{(m)}\right] .
\end{equation*}

(vi) The adapted components of the electromagnetic d-form $F_{\mathcal{P}}^{%
\text{ODEs}}$ produced by the jet first order ODEs system (\ref{ODEs}) and
the pair of Riemannian metrics $\mathcal{P}$ verify the \textbf{generalized
Maxwell equations}%
\begin{equation*}
\left\{ 
\begin{array}{l}
F_{(i)j//1}^{(1)}=\dfrac{1}{4}\mathcal{A}_{\{i,j\}}\left\{ h^{11}\varphi
_{im}\left[ X_{(1)||j//1}^{(m)}-\varphi ^{mr}X_{(1)||r//1}^{(s)}\varphi _{sj}%
\right] \right\} \medskip \\ 
\sum_{\{i,j,k\}}F_{(i)j||k}^{(1)}=0,%
\end{array}%
\right.
\end{equation*}%
where $\mathcal{A}_{\{i,j\}}$ represents an alternate sum, $\sum_{\{i,j,k\}}$
means a cyclic sum and%
\begin{equation*}
F_{(i)j//1}^{(1)}=\dfrac{\partial F_{(i)j}^{(1)}}{\partial t}%
+F_{(i)j}^{(1)}H_{11}^{1}\text{ and }F_{(i)j||k}^{(1)}=\frac{\partial
F_{(i)j}^{(1)}}{\partial x^{k}}-F_{(m)j}^{(1)}\gamma
_{ik}^{m}-F_{(i)m}^{(1)}\gamma _{jk}^{m}
\end{equation*}%
have the geometrical meaning of the horizontal local covariant derivatives $%
"_{//1}"$ and $"_{||k}"$ produced by the Berwald linear connection $B\Gamma
_{0}$ on $J^{1}(T,\mathbb{R}^{n}).$ For more details, please consult [4].

(vii) The \textbf{geometric jet Yang-Mills energy produced by the jet first
order ODEs system (\ref{ODEs}) and the pair of Riemannian metrics} $\mathcal{%
P}$ is defined by the formula%
\begin{equation*}
EYM_{\mathcal{P}}^{\text{ODEs}}(t,x)=\sum_{i=1}^{n-1}\sum_{j=i+1}^{n}\left[
F_{(i)j}^{(1)}\right] ^{2}.
\end{equation*}
\end{theorem}

Now, let us consider on $T\times \mathbb{R}^{n}$ the particular pair of
Euclidian metrics%
\begin{equation*}
\Delta =(h_{11}(t)=1,\varphi _{ij}(x)=\delta _{ij}),
\end{equation*}%
where $\delta _{ij}$ are the Kronecker symbols. Then we obtain the
particular jet least squares Lagrangian function 
\begin{equation*}
JLS_{\Delta }^{\text{ODEs}}:J^{1}(T,\mathbb{R}^{n})\rightarrow \mathbb{R}%
_{+},
\end{equation*}%
defined by%
\begin{eqnarray*}
JLS_{\Delta }^{\text{ODEs}}(t,x^{k},x_{1}^{k}) &=&\delta _{ij}\left[
x_{1}^{i}-X_{(1)}^{(i)}(t,x)\right] \left[ x_{1}^{j}-X_{(1)}^{(j)}(t,x)%
\right] = \\
&=&\sum_{i=1}^{n}\left[ x_{1}^{i}-X_{(1)}^{(i)}(t,x)\right] ^{2}.
\end{eqnarray*}

In this new context, we introduce the following concept:

\begin{definition}
Any geometrical object on $J^{1}(T,\mathbb{R}^{n})$, which is produced by
the jet least squares Lagrangian function $JLS_{\Delta }^{\text{ODEs}}$, via
its attached second order Euler-Lagrange equations, is called \textbf{%
geometrical object produced by the jet first order ODEs system (\ref{ODEs})}.
\end{definition}

As a consequence, particularizing the Theorem \ref{MainThODEs} for the pair
of Euclidian metrics $\mathcal{P}=\Delta $ and taking into account that we
have $H_{11}^{1}(t)=0$ and $\gamma _{ij}^{k}(x)=0$, we immediately get the
following jet geometrical result:

\begin{corollary}
\label{MainCor} (i) The \textbf{canonical non-linear connection on }$J^{1}(T,%
\mathbb{R}^{n})$\textbf{\ produced by the jet first order ODEs system (\ref%
{ODEs})} has the local components%
\begin{equation*}
\Gamma _{\Delta }^{\text{ODEs}}=\left( \bar{M}_{(1)1}^{(i)},\bar{N}%
_{(1)j}^{(i)}\right) ,
\end{equation*}%
where%
\begin{equation*}
\bar{M}_{(1)1}^{(i)}=0\text{ and }\bar{N}_{(1)j}^{(i)}=-\frac{1}{2}\left[ 
\frac{\partial X_{(1)}^{(i)}}{\partial x^{j}}-\frac{\partial X_{(1)}^{(j)}}{%
\partial x^{i}}\right] ,\text{ }\forall \text{ }i,j=\overline{1,n}.
\end{equation*}

(ii) All adapted components of the \textbf{canonical generalized Cartan
connection }$C\Gamma _{\Delta }^{\text{ODEs}}$\textbf{\ produced by the jet
first order ODEs system (\ref{ODEs})} vanish.

(iii) The effective adapted components of the \textbf{torsion} d-tensor 
\textbf{T}$_{\Delta }^{\text{ODEs}}$ of the canonical generalized Cartan
connection $C\Gamma _{\Delta }^{\text{ODEs}}$ \textbf{produced by the jet
first order ODEs system (\ref{ODEs})} are%
\begin{equation*}
\bar{R}_{(1)1j}^{(i)}=\frac{1}{2}\left[ \frac{\partial ^{2}X_{(1)}^{(i)}}{%
\partial t\partial x^{j}}-\frac{\partial ^{2}X_{(1)}^{(j)}}{\partial
t\partial x^{i}}\right] ,\text{ }\forall \text{ }i,j=\overline{1,n},
\end{equation*}%
and%
\begin{equation*}
\bar{R}_{(1)jk}^{(i)}=-\frac{1}{2}\left[ \frac{\partial ^{2}X_{(1)}^{(i)}}{%
\partial x^{k}\partial x^{j}}-\frac{\partial ^{2}X_{(1)}^{(j)}}{\partial
x^{k}\partial x^{i}}\right] ,\text{ }\forall \text{ }i,j,k=\overline{1,n}.
\end{equation*}

(iv) All adapted components of the \textbf{curvature} d-tensor \textbf{R}$%
_{\Delta }^{\text{ODEs}}$ of the canonical generalized Cartan connection $%
C\Gamma _{\Delta }^{\text{ODEs}}$ \textbf{produced by the jet first order
DEs system (\ref{ODEs})} vanish.

(v) The \textbf{geometric electromagnetic distinguished 2-form produced by
the jet first order ODEs system (\ref{ODEs})} has the form%
\begin{equation*}
F_{\Delta }^{\text{ODEs}}=\bar{F}_{(i)j}^{(1)}\delta x_{1}^{i}\wedge dx^{j},
\end{equation*}%
where%
\begin{equation*}
\delta x_{1}^{i}=dx_{1}^{i}+\bar{N}_{(1)k}^{(i)}dx^{k},\text{ }\forall \text{
}i=\overline{1,n},
\end{equation*}%
and%
\begin{equation*}
\bar{F}_{(i)j}^{(1)}=\frac{1}{2}\left[ \frac{\partial X_{(1)}^{(i)}}{%
\partial x^{j}}-\frac{\partial X_{(1)}^{(j)}}{\partial x^{i}}\right] ,\text{ 
}\forall \text{ }i,j=\overline{1,n}.
\end{equation*}

(vi) The adapted components $\bar{F}_{(i)j}^{(1)}$ of the electromagnetic
d-form $F_{\Delta }^{\text{ODEs}}$ produced by the jet first order ODEs
system (\ref{ODEs})\textbf{\ }verify the \textbf{generalized Maxwell
equations}%
\begin{equation*}
\left\{ 
\begin{array}{l}
\bar{F}_{(i)j//1}^{(1)}=\dfrac{1}{4}\mathcal{A}_{\{i,j\}}\left[ \dfrac{%
\partial ^{2}X_{(1)}^{(i)}}{\partial t\partial x^{j}}-\dfrac{\partial
^{2}X_{(1)}^{(j)}}{\partial t\partial x^{i}}\right] =\dfrac{1}{2}\left[ 
\dfrac{\partial ^{2}X_{(1)}^{(i)}}{\partial t\partial x^{j}}-\dfrac{\partial
^{2}X_{(1)}^{(j)}}{\partial t\partial x^{i}}\right] \medskip \\ 
\sum_{\{i,j,k\}}\bar{F}_{(i)j||k}^{(1)}=0,%
\end{array}%
\right.
\end{equation*}%
where $\mathcal{A}_{\{i,j\}}$ represents an alternate sum, $\sum_{\{i,j,k\}}$
means a cyclic sum and%
\begin{equation*}
\bar{F}_{(i)j//1}^{(1)}=\dfrac{\partial \bar{F}_{(i)j}^{(1)}}{\partial t}%
\text{ and }\bar{F}_{(i)j||k}^{(1)}=\frac{\partial \bar{F}_{(i)j}^{(1)}}{%
\partial x^{k}},\text{ }\forall \text{ }i,j,k=\overline{1,n}.
\end{equation*}

(vii) The \textbf{geometric jet Yang-Mills energy produced by the jet first
order ODEs system (\ref{ODEs})} has the expression%
\begin{equation*}
EYM_{\Delta }^{\text{ODEs}}(t,x)=\sum_{i=1}^{n-1}\sum_{j=i+1}^{n}\left[ \bar{%
F}_{(i)j}^{(1)}\right] ^{2}=\frac{1}{4}\sum_{i=1}^{n-1}\sum_{j=i+1}^{n}\left[
\frac{\partial X_{(1)}^{(i)}}{\partial x^{j}}-\frac{\partial X_{(1)}^{(j)}}{%
\partial x^{i}}\right] ^{2}.
\end{equation*}
\end{corollary}

\begin{remark}
\label{Formulas} If we use the matriceal notations

\begin{itemize}
\item $J\left( X_{(1)}\right) =\left( \dfrac{\partial X_{(1)}^{(i)}}{%
\partial x^{j}}\right) _{i,j=\overline{1,n}}$ - the \textbf{Jacobian matrix},

\item $\bar{N}_{(1)}=\left( \bar{N}_{(1)j}^{(i)}\right) _{i,j=\overline{1,n}%
} $ - the \textbf{non-linear connection matrix},

\item $\bar{R}_{(1)1}=\left( \bar{R}_{(1)1j}^{(i)}\right) _{i,j=\overline{1,n%
}},$ - the \textbf{temporal torsion matrix},

\item $\bar{R}_{(1)k}=\left( \bar{R}_{(1)jk}^{(i)}\right) _{i,j=\overline{1,n%
}},$ $\forall $ $k=\overline{1,n},$ - the \textbf{spatial torsion matrices},

\item $\bar{F}^{(1)}=\left( \bar{F}_{(i)j}^{(1)}\right) _{i,j=\overline{1,n}%
} $ - the \textbf{electromagnetic matrix},
\end{itemize}

then the following matriceal geometrical relations attached to the jet first
order ODEs system (\ref{ODEs}) hold good:

\begin{enumerate}
\item $\bar{N}_{(1)}=-\dfrac{1}{2}\left[ J\left( X_{(1)}\right) -\text{ }%
^{T}J\left( X_{(1)}\right) \right] ;$

\item $\bar{R}_{(1)1}=-\dfrac{\partial }{\partial t}\left[ \bar{N}_{(1)}%
\right] ;$

\item $\bar{R}_{(1)k}=\dfrac{\partial }{\partial x^{k}}\left[ \bar{N}_{(1)}%
\right] ,$ $\forall $ $k=\overline{1,n};$

\item $\bar{F}^{(1)}=-\bar{N}_{(1)};$

\item $EYM_{\Delta }^{\text{ODEs}}(t,x)=\dfrac{1}{2}\cdot Trace\left[ \bar{F}%
^{(1)}\cdot \text{ }^{T}\bar{F}^{(1)}\right] ,$

that is the jet electromagnetic Yang-Mills energy coincides with the square
of the norm of the skew-symmetric electromagnetic matrix $\bar{F}^{(1)}$ in
the Lie algebra $o(n)=L(O(n)).$
\end{enumerate}
\end{remark}

\begin{remark}
Note that the spatial torsion matrix $\bar{R}_{(1)k}$ does not coincide for $%
k=1$ with the temporal torsion matrix $\bar{R}_{(1)1}$. We have only an
overlap of notations.
\end{remark}

\section{Jet Riemann-Lagrange geometry produced by a non-homogenous linear
ODEs system of order one}

\hspace{4mm} In this Section we apply the preceding jet Riemann-Lagrange
geometrical results for a non-homogenous linear ODEs system of order one. In
this way, let us consider the following non-homogenous linear first order
ODEs system locally described, in a convenient chart on $J^{1}(T,\mathbb{R}%
^{n})$, by the differential equations%
\begin{equation}
\frac{dx^{i}}{dt}=\sum_{k=1}^{n}a_{(1)k}^{(i)}(t)x^{k}+f_{(1)}^{(i)}(t),%
\text{ }\forall \text{ }i=\overline{1,n},  \label{LODEs}
\end{equation}%
where the local components $a_{(1)k}^{(i)}$ and $f_{(1)}^{(i)}$ transform
after the tensorial rules%
\begin{equation*}
a_{(1)k}^{(i)}=\dfrac{\partial x^{i}}{\partial \widetilde{x}^{j}}\dfrac{d%
\widetilde{t}}{dt}\cdot \widetilde{a}_{(1)k}^{(j)},\text{ }\forall \text{ }k=%
\overline{1,n},
\end{equation*}%
and%
\begin{equation*}
f_{(1)}^{(i)}=\frac{\partial x^{i}}{\partial \widetilde{x}^{j}}\frac{d%
\widetilde{t}}{dt}\cdot \widetilde{f}_{(1)}^{(j)}.
\end{equation*}

\begin{remark}
We suppose that the product manifold $T\times \mathbb{R}^{n}\subset J^{1}(T,%
\mathbb{R}^{n})$ is endowed \textbf{a priori} with the pair of Euclidian
metrics $\Delta =(1,\delta _{ij}),$ with respect to the coordinates $%
(t,x^{i})$.
\end{remark}

It is obvious that the non-homogenous linear ODEs system (\ref{LODEs}) is a
particular case of the jet first order non-linear ODEs system (\ref{ODEs})
for%
\begin{equation}
X_{(1)}^{(i)}(t,x)=\sum_{k=1}^{n}a_{(1)k}^{(i)}(t)x^{k}+f_{(1)}^{(i)}(t),%
\text{ }\forall \text{ }i=\overline{1,n}.  \label{XLODEs}
\end{equation}

In order to expose the main jet Riemann-Lagrange geometrical objects that
characterize the non-homogenous linear ODEs system (\ref{LODEs}), we use the
matriceal notation%
\begin{equation*}
A_{(1)}=\left( a_{(1)j}^{(i)}(t)\right) _{i,j=\overline{1,n}}.
\end{equation*}

In this context, applying our preceding jet geometrical Riemann-Lagrange
theory to the non-homogenous linear ODEs system (\ref{LODEs}) and the pair
of Euclidian metrics $\Delta =(1,\delta _{ij})$, we get:

\begin{theorem}
(i) The \textbf{canonical non-linear connection on }$J^{1}(T,\mathbb{R}^{n})$%
\textbf{\ produced by the non-homogenous linear ODEs system (\ref{LODEs})}
has the local components%
\begin{equation*}
\hat{\Gamma}=\left( 0,\hat{N}_{(1)j}^{(i)}\right) ,
\end{equation*}%
where $\hat{N}_{(1)j}^{(i)}$ are the entries of the matrix 
\begin{equation*}
\hat{N}_{(1)}=\left( \hat{N}_{(1)j}^{(i)}\right) _{i,j=\overline{1,n}}=-%
\frac{1}{2}\left[ A_{(1)}-\text{ }^{T}A_{(1)}\right] .
\end{equation*}

(ii) All adapted components of the \textbf{canonical generalized Cartan
connection }$C\hat{\Gamma}$\textbf{\ produced by the non-homogenous linear
ODEs system (\ref{LODEs})} vanish.

(iii) The effective adapted components $\hat{R}_{(1)1j}^{(i)}$ of the 
\textbf{torsion} d-tensor \textbf{\^{T}} of the canonical generalized Cartan
connection $C\hat{\Gamma}$ \textbf{produced by the non-homogenous linear
ODEs system (\ref{LODEs}) }are the entries of the matrices%
\begin{equation*}
\hat{R}_{(1)1}=\left( \hat{R}_{(1)1j}^{(i)}\right) _{i,j=\overline{1,n}}=%
\frac{1}{2}\left[ \dot{A}_{(1)}-\text{ }^{T}\dot{A}_{(1)}\right] ,
\end{equation*}%
where%
\begin{equation*}
\dot{A}_{(1)}=\frac{d}{dt}\left[ A_{(1)}\right] .
\end{equation*}

(iv) All adapted components of the \textbf{curvature} d-tensor \textbf{\^{R}}
of the canonical generalized Cartan connection $C\hat{\Gamma}$ \textbf{%
produced by the non-homogenous linear ODEs system (\ref{LODEs})} vanish.

(v) The \textbf{geometric electromagnetic distinguished 2-form produced by
the non-homogenous linear ODEs system (\ref{LODEs})} is given by%
\begin{equation*}
\hat{F}=\hat{F}_{(i)j}^{(1)}\delta x_{1}^{i}\wedge dx^{j},
\end{equation*}%
where%
\begin{equation*}
\delta x_{1}^{i}=dx_{1}^{i}-\frac{1}{2}\left[ a_{(1)k}^{(i)}-a_{(1)i}^{(k)}%
\right] dx^{k},\text{ }\forall \text{ }i=\overline{1,n},
\end{equation*}%
and the adapted components $\hat{F}_{(i)j}^{(1)}$ are the entries of the
matrix%
\begin{equation*}
\hat{F}^{(1)}=\left( \hat{F}_{(i)j}^{(1)}\right) _{i,j=\overline{1,n}}=-\hat{%
N}_{(1)}=\frac{1}{2}\left[ A_{(1)}-\text{ }^{T}A_{(1)}\right] ,
\end{equation*}%
that is 
\begin{equation*}
\hat{F}_{(i)j}^{(1)}=\frac{1}{2}\left[ a_{(1)j}^{(i)}-a_{(1)i}^{(j)}\right] .
\end{equation*}

(vi) The \textbf{\ jet Yang-Mills energy produced by the non-homogenous
linear ODEs system (\ref{LODEs})} is given by the formula%
\begin{equation*}
EYM^{\text{NHLODEs}}(t)=\frac{1}{4}\sum_{i=1}^{n-1}\sum_{j=i+1}^{n}\left[
a_{(1)j}^{(i)}-a_{(1)i}^{(j)}\right] ^{2}.
\end{equation*}
\end{theorem}

\begin{proof}
Using the relations (\ref{XLODEs}), we easily deduce that we have the
Jacobian matrix%
\begin{equation*}
J\left( X_{(1)}\right) =A_{(1)}.
\end{equation*}

Consequently, applying the Corollary \ref{MainCor} to the non-homogenous
linear ODEs system (\ref{LODEs}), together with the Remark \ref{Formulas},
we obtain the required results.
\end{proof}

\begin{remark}
The entire jet Riemann-Lagrange geometry produced by the non-homogenous
linear ODEs system (\ref{LODEs}) does not depend on the non-homogeneity
terms $f_{(1)}^{(i)}(t)$.
\end{remark}

\begin{remark}
The \textbf{\ }jet Yang-Mills energy produced by the non-homogenous linear
ODEs system (\ref{LODEs}) vanishes if and only if the matrix $A_{(1)}$ is a
symmetric one. In this case, the entire jet Riemann-Lagrange geometry
produced by the non-homogenous linear ODEs system (\ref{LODEs}) vanish, so
it does not offer geometrical informations about the system (\ref{LODEs}).
However, it is important to note that in this particular situation we have
the symetry of the matrix $A_{(1)},$ which implies that the matrix $A_{(1)}$
is diagonalizable.
\end{remark}

\begin{remark}
All torsion adapted components of a non-homogenous linear ODEs system with
constant coefficients $a_{(1)j}^{(i)}$ are zero.
\end{remark}

\section{Jet Riemann-Lagrange geometry produced by a superior order ODE}

\hspace{4mm} Let us consider the superior order ODE expressed by%
\begin{equation}
y^{(n)}(t)=f(t,y(t),y^{\prime }(t),...,y^{(n-1)}(t)),\text{ }n\geq 2,
\label{SODE}
\end{equation}%
where $y(t)$ is an unknown function, $y^{(k)}(t)$ is the derivative of order 
$k$ of the unknown function $y(t)$ for each $k\in \{0,1,...,n\}$ and $f$ is
a given differentiable function depending on the distinct variables $t,$ $%
y(t),$ $y^{\prime }(t),...,$ $y^{(n-1)}(t).$

It is well known the fact that, using the notations%
\begin{equation*}
x^{1}=y,\text{ }x^{2}=y^{\prime },\text{ }...,\text{ }x^{n}=y^{(n-1)},
\end{equation*}%
the superior order ODE (\ref{SODE}) is equivalent with the non-linear ODEs
system of order one%
\begin{equation}
\left\{ 
\begin{array}{l}
\dfrac{dx^{1}}{dt}=x^{2}\medskip \\ 
\dfrac{dx^{2}}{dt}=x^{3}\medskip \\ 
\cdot \\ 
\cdot \\ 
\cdot \\ 
\dfrac{dx^{n-1}}{dt}=x^{n}\medskip \\ 
\dfrac{dx^{n}}{dt}=f(t,x^{1},x^{2},...,x^{n}).%
\end{array}%
\right.  \label{NLSODE}
\end{equation}

But, the first order non-linear ODEs system (\ref{NLSODE}) can be regarded,
in a convenient local chart, as a particular case of the jet non-linear ODEs
system of order one (\ref{ODEs}), taking%
\begin{equation}
\begin{array}{lll}
X_{(1)}^{(1)}(t,x)=x^{2}, & X_{(1)}^{(2)}(t,x)=x^{3}, & \cdot \cdot \cdot
\medskip \\ 
\cdot \cdot \cdot & X_{(1)}^{(n-1)}(t,x)=x^{n}, & 
X_{(1)}^{(n)}(t,x)=f(t,x^{1},x^{2},...,x^{n}),%
\end{array}
\label{XSODE}
\end{equation}%
where we suppose that the geometrical object $X=\left(
X_{(1)}^{(i)}(t,x)\right) $ behaves as a d-tensor on $J^{1}(T,\mathbb{R}%
^{n}) $.

\begin{remark}
We assume that the product manifold $T\times \mathbb{R}^{n}\subset J^{1}(T,%
\mathbb{R}^{n})$ is endowed \textbf{a priori} with the pair of Euclidian
metrics $\Delta =(1,\delta _{ij}),$ with respect to the coordinates $%
(t,x^{i})$.
\end{remark}

\begin{definition}
Any geometrical object on $J^{1}(T,\mathbb{R}^{n})$, which is produced by
the \textbf{\ first order non-linear ODEs system (\ref{NLSODE})} is called 
\textbf{geometrical object produced by the superior order ODE (\ref{SODE})}.
\end{definition}

In this context, the Riemann-Lagrange geometrical behavior on the 1-jet
space $J^{1}(T,\mathbb{R}^{n})$ of the superior order ODE (\ref{SODE}) is
described in the following result:

\begin{theorem}
\label{MainThSODE} (i) The \textbf{canonical non-linear connection on }$%
J^{1}(T,\mathbb{R}^{n})$\textbf{\ produced by the superior order ODE (\ref%
{SODE})} has the local components%
\begin{equation*}
\check{\Gamma}=\left( 0,\check{N}_{(1)j}^{(i)}\right) ,
\end{equation*}%
where $\check{N}_{(1)j}^{(i)}$ are the entries of the matrix $\check{N}%
_{(1)}=\left( \check{N}_{(1)j}^{(i)}\right) _{i,j=\overline{1,n}}=$%
\begin{equation*}
=-\frac{1}{2}\left( 
\begin{array}{cccccccc}
0 & 1 & 0 & \cdot & \cdot & 0 & 0 & -\dfrac{\partial f}{\partial x^{1}}%
\medskip \\ 
-1 & 0 & 1 & \cdot & \cdot & 0 & 0 & -\dfrac{\partial f}{\partial x^{2}}%
\medskip \\ 
0 & -1 & 0 & \cdot & \cdot & 0 & 0 & -\dfrac{\partial f}{\partial x^{3}}%
\medskip \\ 
\cdot & \cdot & \cdot & \cdot & \cdot & \cdot & \cdot & \cdot \medskip \\ 
\cdot & \cdot & \cdot & \cdot & \cdot & \cdot & \cdot & \cdot \medskip \\ 
0 & 0 & 0 & \cdot & \cdot & 0 & 1 & -\dfrac{\partial f}{\partial x^{n-2}}%
\medskip \\ 
0 & 0 & 0 & \cdot & \cdot & -1 & 0 & 1-\dfrac{\partial f}{\partial x^{n-1}}%
\medskip \\ 
\dfrac{\partial f}{\partial x^{1}} & \dfrac{\partial f}{\partial x^{2}} & 
\dfrac{\partial f}{\partial x^{3}} & \cdot & \cdot & \dfrac{\partial f}{%
\partial x^{n-2}} & -1+\dfrac{\partial f}{\partial x^{n-1}} & 0%
\end{array}%
\right) .
\end{equation*}

(ii) All adapted components of the \textbf{canonical generalized Cartan
connection }$C\check{\Gamma}$\textbf{\ produced by the superior order ODE (%
\ref{SODE})} vanish.

(iii) The effective adapted components of the \textbf{torsion} d-tensor 
\textbf{\v{T}} of the canonical generalized Cartan connection $C\check{\Gamma%
}$ \textbf{produced by the superior order ODE (\ref{SODE})} are the entries
of the matrices%
\begin{equation*}
\check{R}_{(1)1}=\frac{1}{2}\left( 
\begin{array}{cccccc}
0 & 0 & \cdot & \cdot & 0 & -\dfrac{\partial ^{2}f}{\partial t\partial x^{1}}%
\medskip \\ 
0 & 0 & \cdot & \cdot & 0 & -\dfrac{\partial ^{2}f}{\partial t\partial x^{2}}%
\medskip \\ 
\cdot & \cdot & \cdot & \cdot & \cdot & \cdot \medskip \\ 
\cdot & \cdot & \cdot & \cdot & \cdot & \cdot \medskip \\ 
0 & 0 & \cdot & \cdot & 0 & -\dfrac{\partial ^{2}f}{\partial t\partial
x^{n-1}}\medskip \\ 
\dfrac{\partial ^{2}f}{\partial t\partial x^{1}} & \dfrac{\partial ^{2}f}{%
\partial t\partial x^{2}} & \cdot & \cdot & \dfrac{\partial ^{2}f}{\partial
t\partial x^{n-1}} & 0%
\end{array}%
\right)
\end{equation*}%
and%
\begin{equation*}
\check{R}_{(1)k}=-\frac{1}{2}\left( 
\begin{array}{cccccc}
0 & 0 & \cdot & \cdot & 0 & -\dfrac{\partial ^{2}f}{\partial x^{k}\partial
x^{1}}\medskip \\ 
0 & 0 & \cdot & \cdot & 0 & -\dfrac{\partial ^{2}f}{\partial x^{k}\partial
x^{2}}\medskip \\ 
\cdot & \cdot & \cdot & \cdot & \cdot & \cdot \medskip \\ 
\cdot & \cdot & \cdot & \cdot & \cdot & \cdot \medskip \\ 
0 & 0 & \cdot & \cdot & 0 & -\dfrac{\partial ^{2}f}{\partial x^{k}\partial
x^{n-1}}\medskip \\ 
\dfrac{\partial ^{2}f}{\partial x^{k}\partial x^{1}} & \dfrac{\partial ^{2}f%
}{\partial x^{k}\partial x^{2}} & \cdot & \cdot & \dfrac{\partial ^{2}f}{%
\partial x^{k}\partial x^{n-1}} & 0%
\end{array}%
\right) ,
\end{equation*}%
where $k\in \{1,2,...,n\}$.

(iv) All adapted components of the \textbf{curvature} d-tensor \textbf{\v{R}}
of the canonical generalized Cartan connection $C\check{\Gamma}$ \textbf{%
produced by the superior order ODE (\ref{SODE})} vanish.

(v) The \textbf{geometric electromagnetic distinguished 2-form produced by
the superior order ODE (\ref{SODE})} has the form%
\begin{equation*}
\check{F}=\check{F}_{(i)j}^{(1)}\delta x_{1}^{i}\wedge dx^{j},
\end{equation*}%
where%
\begin{equation*}
\delta x_{1}^{i}=dx_{1}^{i}+\check{N}_{(1)k}^{(i)}dx^{k},\text{ }\forall 
\text{ }i=\overline{1,n},
\end{equation*}%
and the adapted components $\check{F}_{(i)j}^{(1)}$ are the entries of the
matrix%
\begin{equation*}
\check{F}^{(1)}=\left( \check{F}_{(i)j}^{(1)}\right) _{i,j=\overline{1,n}}=-%
\check{N}_{(1)}.
\end{equation*}

(vi) The \textbf{jet geometric Yang-Mills energy produced by the superior
order ODE (\ref{SODE})} is given by the formula%
\begin{equation*}
EYM^{\text{SODE}}(t,x)=\frac{1}{4}\left[ n-1-2\frac{\partial f}{\partial
x^{n-1}}+\sum_{j=1}^{n-1}\left( \frac{\partial f}{\partial x^{j}}\right) ^{2}%
\right] .
\end{equation*}
\end{theorem}

\begin{proof}
By partial derivatives, the relations (\ref{XSODE}) lead to the Jacobian
matrix%
\begin{equation*}
J\left( X_{(1)}\right) =\left( 
\begin{array}{ccccccc}
0 & 1 & 0 & \cdot & \cdot & 0 & 0\medskip \\ 
0 & 0 & 1 & \cdot & \cdot & 0 & 0\medskip \\ 
\cdot & \cdot & \cdot & \cdot & \cdot & \cdot & \cdot \medskip \\ 
\cdot & \cdot & \cdot & \cdot & \cdot & \cdot & \cdot \medskip \\ 
0 & 0 & 0 & \cdot & \cdot & 0 & 1\medskip \\ 
\dfrac{\partial f}{\partial x^{1}} & \dfrac{\partial f}{\partial x^{2}} & 
\dfrac{\partial f}{\partial x^{3}} & \cdot & \cdot & \dfrac{\partial f}{%
\partial x^{n-1}} & \dfrac{\partial f}{\partial x^{n}}%
\end{array}%
\right) .
\end{equation*}

In conclusion, the Corollary \ref{MainCor}, together with the Remark \ref%
{Formulas}, applied to first order non-linear ODEs system (\ref{NLSODE}),%
\textbf{\ }give what we were looking for.
\end{proof}

\section{Riemann-Lagrange geometry produced by a non-homogenous linear ODE
of superior order}

\hspace{4mm} If we consider the non-homogenous linear ODE of order $n\in 
\mathbb{N}$, $n\geq 2$, expressed by%
\begin{equation}
a_{0}(t)y^{(n)}+a_{1}(t)y^{(n-1)}+...+a_{n-1}(t)y^{\prime }+a_{n}(t)y=b(t),
\label{NHLSODE}
\end{equation}%
where $b(t)$ and $a_{i}(t)$, $\forall $ $i=\overline{0,n}$, are given
differentiable real functions and $a_{0}(t)\neq 0$, $\forall $ $t\in \lbrack
a,b]$, then we recover the superior order ODE (\ref{SODE}) for the
particular function%
\begin{equation}
f(t,x)=\frac{b(t)}{a_{0}(t)}-\frac{a_{n}(t)}{a_{0}(t)}\cdot x^{1}-\frac{%
a_{n-1}(t)}{a_{0}(t)}\cdot x^{2}-...-\frac{a_{1}(t)}{a_{0}(t)}\cdot x^{n},
\label{fNHLSODE}
\end{equation}%
where we recall that we have%
\begin{equation*}
y=x^{1},y^{\prime }=x^{2},...,y^{(n-1)}=x^{n}.
\end{equation*}

Consequently, we can derive the jet Riemann-Lagrange geometry attached to
the non-homogenous linear superior order ODE (\ref{NHLSODE}).

\begin{corollary}
(i) The \textbf{canonical non-linear connection on }$J^{1}(T,\mathbb{R}^{n})$%
\textbf{\ produced by the non-homogenous linear superior order ODE (\ref%
{NHLSODE})} has the local components%
\begin{equation*}
\tilde{\Gamma}=\left( 0,\tilde{N}_{(1)j}^{(i)}\right) ,
\end{equation*}%
where $\tilde{N}_{(1)j}^{(i)}$ are the entries of the matrix 
\begin{eqnarray*}
\tilde{N}_{(1)} &=&\left( \tilde{N}_{(1)j}^{(i)}\right) _{i,j=\overline{1,n}%
}= \\
&=&-\frac{1}{2}\left( 
\begin{array}{cccccccc}
0 & 1 & 0 & \cdot  & \cdot  & 0 & 0 & \dfrac{a_{n}}{a_{0}}\medskip  \\ 
-1 & 0 & 1 & \cdot  & \cdot  & 0 & 0 & \dfrac{a_{n-1}}{a_{0}}\medskip  \\ 
0 & -1 & 0 & \cdot  & \cdot  & 0 & 0 & \dfrac{a_{n-2}}{a_{0}}\medskip  \\ 
\cdot  & \cdot  & \cdot  & \cdot  & \cdot  & \cdot  & \cdot  & \cdot
\medskip  \\ 
\cdot  & \cdot  & \cdot  & \cdot  & \cdot  & \cdot  & \cdot  & \cdot
\medskip  \\ 
0 & 0 & 0 & \cdot  & \cdot  & 0 & 1 & \dfrac{a_{3}}{a_{0}}\medskip  \\ 
0 & 0 & 0 & \cdot  & \cdot  & -1 & 0 & 1+\dfrac{a_{2}}{a_{0}}\medskip  \\ 
-\dfrac{a_{n}}{a_{0}} & -\dfrac{a_{n-1}}{a_{0}} & -\dfrac{a_{n-2}}{a_{0}} & 
\cdot  & \cdot  & -\dfrac{a_{3}}{a_{0}} & -1-\dfrac{a_{2}}{a_{0}} & 0%
\end{array}%
\right) .
\end{eqnarray*}

(ii) All adapted components of the \textbf{canonical generalized Cartan
connection }$C\tilde{\Gamma}$\textbf{\ produced by the non-homogenous linear
superior order ODE (\ref{NHLSODE})} vanish.

(iii) All adapted components of the \textbf{torsion} d-tensor \textbf{\~{T}}
of the canonical generalized Cartan connection $C\tilde{\Gamma}$ \textbf{%
produced by the non-homogenous linear superior order ODE (\ref{NHLSODE})}
are zero, except the temporal components%
\begin{equation*}
\tilde{R}_{(1)1n}^{(i)}=-\tilde{R}_{(1)1i}^{(n)}=\frac{a_{n-i+1}^{\prime
}a_{0}-a_{n-i+1}a_{0}^{\prime }}{2a_{0}^{2}},\text{ }\forall \text{ }i=%
\overline{1,n-1},
\end{equation*}%
where we denoted by $"$ $^{\prime }$ $"$ the derivatives of the functions $%
a_{k}(t)$.

(iv) All adapted components of the \textbf{curvature} d-tensor \textbf{\~{R}}
of the canonical generalized Cartan connection $C\tilde{\Gamma}$ \textbf{%
produced by the non-homogenous linear superior order ODE (\ref{NHLSODE})}
vanish.

(v) The \textbf{geometric electromagnetic distinguished 2-form produced by
the non-homogenous linear superior order ODE (\ref{NHLSODE})} has the
expression%
\begin{equation*}
\tilde{F}=\tilde{F}_{(i)j}^{(1)}\delta x_{1}^{i}\wedge dx^{j},
\end{equation*}%
where%
\begin{equation*}
\delta x_{1}^{i}=dx_{1}^{i}+\tilde{N}_{(1)k}^{(i)}dx^{k},\text{ }\forall 
\text{ }i=\overline{1,n},
\end{equation*}%
and the adapted components $\tilde{F}_{(i)j}^{(1)}$ are the entries of the
matrix%
\begin{equation*}
\tilde{F}^{(1)}=\left( \tilde{F}_{(i)j}^{(1)}\right) _{i,j=\overline{1,n}}=-%
\tilde{N}_{(1)}.
\end{equation*}

(vi) The \textbf{jet geometric Yang-Mills electromagnetic energy produced by
the non-homogenous linear superior order ODE (\ref{NHLSODE})} has the form%
\begin{equation*}
EYM^{\text{NHLSODE}}(t)=\frac{1}{4}\left[ n-1+2\frac{a_{2}}{a_{0}}%
+\sum_{j=2}^{n}\frac{a_{j}^{2}}{a_{0}^{2}}\right] .
\end{equation*}
\end{corollary}

\begin{proof}
We apply the Theorem \ref{MainThSODE} for the particular function (\ref%
{fNHLSODE}) and we use the relations%
\begin{equation*}
\frac{\partial f}{\partial x^{j}}=-\frac{a_{n-j+1}}{a_{0}},\text{ }\forall 
\text{ }j=\overline{1,n}.
\end{equation*}
\end{proof}

\begin{remark}
The entire jet Riemann-Lagrange geometry produced by the non-homogenous
linear superior order ODE (\ref{NHLSODE}) is independent by the term of
non-homogeneity $b(t)$. In author's opinion, this fact emphasizes that the
most important role in the study of the ODE (\ref{NHLSODE}) is played by its
attached homogenous linear superior order ODE.
\end{remark}

\begin{example}
The law of motion without friction (\textbf{harmonic oscillator}) of a
material point of mass $m>0$, which is placed on a spring having the
constant of elasticity $k>0$, is given by the homogenous linear ODE of order
two%
\begin{equation}
\frac{d^{2}y}{dt^{2}}+\omega ^{2}y=0,  \label{Oscilator}
\end{equation}%
where the coordinate $y$ measures the distance from the mass's equlibrium
point and $\omega ^{2}=k/m.$ It follows that we have%
\begin{equation*}
n=2,\text{ }a_{0}(t)=1,\text{ }a_{1}(t)=0\text{ and }a_{2}(t)=\omega ^{2},
\end{equation*}%
that is the \textbf{harmonic oscillator} second order ODE (\ref{Oscilator})
provides the \textbf{jet geometric Yang-Mills electromagnetic energy}%
\begin{equation*}
EYM^{\text{Harmonic Oscillator}}=\frac{1}{4}\left( 1+\omega ^{2}\right) ^{2}.
\end{equation*}
\end{example}

\textbf{Open problem.} There exists a real physical interpretation for the
previous jet geometric Yang-Mills electromagnetic energy attached to the
harmonic oscillator?

\textbf{Acknowledgements.} A version of this paper was presented at \textit{%
Conference of Differential Geometry dedicated to the Memory of Professor
Kostake Teleman (1933-2007)}, University of Bucharest, May 15-17, 2009.

The present work was supported by Contract with Sinoptix No. 8441/2009.

\textbf{Author's address:} Mircea N{\scriptsize EAGU}

University Transilvania of Bra\c{s}ov

Faculty of Mathematics and Informatics

Department of Algebra, Geometry and Differential Equations

B-dul Eroilor, Nr. 29, 500036 Bra\c{s}ov, Romania.

E-mail: mircea.neagu@unitbv.ro

Website: http://www.2collab.com/user:mirceaneagu

\end{document}